\documentclass[reqno,10pt]{amsart}
\usepackage{amsmath,amsthm,amssymb}
\usepackage{setspace}
\usepackage[pdftex]{hyperref}
\theoremstyle{definition}
 \newtheorem{defi}{Definition}[section]

\theoremstyle{plain}
 \newtheorem{prop}{Proposition}[section]

\theoremstyle{plain}
 \newtheorem{thm}{Theorem}[section]

\theoremstyle{remark}
 \newtheorem{exam}{Example}[section]

\theoremstyle{plain}
 \newtheorem{lem}{Lemma}[section]

\theoremstyle{plain}
 \newtheorem{cor}{Corollary}[section]

\theoremstyle{remark}
 
\hypersetup{bookmarksopen,bookmarksnumbered,colorlinks,%
linkcolor=blue,legalpaper,pagebackref,%
citecolor=red,%
urlcolor=blue,%
pdftitle=Some Divisibility Properties in Ring of Polynomials over a Unique Factorization Domain,%
pdfauthor=Luis F. C\'aceres \& Jos\'e V\'elez-Marulanda \textcopyright,pdfsubject=Mathematics}

\newcommand{\Z}{\mathbb{Z}}
\newcommand{\Q}{\mathbb{Q}}

\newcommand{\Norm}{\textit{\textbf{N}}}
\newcommand{\R}{\mathbb{R}}
\newcommand{\Prime}{\mathfrak{P}}

\newcommand{\C}{\mathfrak{C}}
\newcommand{\J}{\mathfrak{J}}

\newcommand{\M}{\mathfrak{M}}
\newcommand{\A}{\mathfrak{A}}

\newcommand{\LA}{\left\langle}
\newcommand{\RA}{\right\rangle}
\newcommand{\Zd}{\mathbb{Z}[\sqrt{d}]}
\newcommand{\Zdn}{\mathbb{Z}[\sqrt{d_1},\ldots,\sqrt{d_n}]}

\DeclareMathOperator{\lcm}{lcm}
\title{Some Divisibility Properties in Ring of Polynomials over a UFD}
\author{Luis F. C\'aceres and Jos\'e A. V\'elez-Marulanda}
\address{Luis F. C\'aceres:}
\address{Department of Mathematics, University of Puerto Rico at Mayag\"uez}
\email{\texttt{lcaceres@math.uprm.edu}}
\urladdr{\texttt{http://www.math.uprm.edu/$\sim$lcaceres/}}
\curraddr{P.O.Box 5622 Mayag\"uez, Puerto Rico 00681-5622}
\address{Jos\'e A. V\'elez-Marulanda:}
\address{Graduate Student, University of Iowa}
\email{\texttt{jose-velezmarulanda@uiowa.edu}}
\curraddr{14 MacLean Hall, Iowa City, Iowa 52242-1419}

\renewcommand{\S}{\mathfrak{S}}
\date{\today}
\begin{document}
\bibliographystyle{plain}
\maketitle 
\begin{abstract}
Using polynomial evaluation, we give some useful criteria to answer questions about divisibility of polynomials. This allows us to develop interesting results concerning the prime elements in the domain of coefficients. In particular, it is possible to prove that under certain conditions, the domain of coefficients must have infinitely many prime elements. We give alternative characterizations for $D-$rings and present various examples.
\end{abstract}
\section{Introduction}
An interesting question about divisibility of polynomials is the following: given
$f(x)$ and $g(x)$ polynomials with coefficients in the ring of integers $\Z$ such that $f(n)|g(n)$ for all
$n\in \Z$, does one have that $f(x)|g(x)$ in $\Z[x]$?
Take for example $f(x)= 5$, and
$g(x)= x^5-x$; by Fermat's Little Theorem
we have that for all $n\in \Z$, $5|n^5-n$ in $\Z$, but clearly $5\nmid x^5-x$ in $\Z[x]$. However, $\Z$ satisfies some
properties showing that in many nontrivial cases the answer to
that question is affirmative. In order to solve this interrogant, we study some
divisibility properties in arbitrary \textit{unique factorization domains} ($UFD$), namely: \textit{infinite
primes property (IPP), degree polynomial property (DPP),
evaluation polynomial property (EPP)} and \textit{strong
evaluation polynomial property (SEPP)}. These properties provide us
useful tools to understand divisibility in the ring
$\Z[x]$ and in any ring of polynomials $D[x]$ for any $UFD$ $D$.
Another property that will be useful is the \textit{$D$-ring} property. In
Section 3 we study this property in detail, we give many examples
and we prove that in a $UFD$, all these properties are equivalent. In the last section we provide some examples.
\section{Basic Definitions}
\begin{defi}\label{IPP}
An integral domain $D$ satisfies the \textit{\textbf{infinite
primes property (IPP)}} if given $g(x)\in D[x]$ with $\deg g(x)
\geq 1$ the set
\[\{ p\in P : (\exists k\in D)(g(k) \not= 0 \text{ and } p \vert g(k)\}\]
is infinite, where $P$ is the set of primes in $D$.
\end{defi}
It is clear that fields do not satisfy $IPP$ (there are no primes in
fields!). It also follows from the definition that rings
satisfying the $IPP$ property must contain infinitely many primes.
\begin{exam}
Let $g(x)=(x-3)(x+2)\in \Z[x]$. Note that $g(3)=0$. Let $p$ be a prime
such that $p|g(p+3)=p(p+5)$. Note that $\Z$ has infinitely many
primes satisfying this condition. Then
\[\{p\in P: (\exists k\in \Z)(g(k)\not=0 \text{ and } p\vert g(k)\},\]
where $P$ is the set of primes of $\Z$, is infinite.
In general, given $g(x)\in \Z[x]$ such that $g(a)=0$ for some
$a\in \Z$, the set
\[\{p\in P: (\exists k\in \Z)(g(k)\not=0 \text{ and } p\vert g(k)\},\]
where $P$ is the set of primes of $\Z$, is infinite. See proof of
Proposition \ref{prop1} below.
\end{exam}
\begin{exam}
Let $p$ be a prime in $\Z$ such that $p\equiv 1\mod 4$. It is well-known (see
\cite[pg~151]{herstein}) that we can find an integer $k$ such that $k^2+1\equiv
0\mod p$. It is also well-known that there are infinitely many primes $p$ such that
$p\equiv 1\mod 4$ (see \cite{burton}). Therefore, the set
\[\{ p\in P : (\exists k\in \Z)(g(k) \not= 0 \text{ and } p \vert g(k)\},\]
where $g(x)=x^2+1$ and $P$ is the set of primes of $\Z$, is
infinite.
\end{exam}
\begin{exam}
Consider the polynomial $g(x)=x^2-2$. The congruence $x^2\equiv 2 \mod p$ has solution if and only if $p\equiv 1\mod 8$. It is well-known that the set of primes of the form $p\equiv 1\mod 8$ is infinite. Hence, the set
\[\{ p\in P : (\exists k\in \Z)(g(k) \not= 0 \text{ and } p \vert g(k)\},\]
where $g(x)=x^2-2$ and $P$ is the set of primes of $\Z$, is
infinite.
\end{exam}
We show that the ring of integers $\Z$ satisfies $IPP$.
\begin{lem}\label{ZsatisfiesIPP}
The ring of integers $\Z$ satisfies $IPP$.
\end{lem}
\begin{proof}
Let $f(x)\in \Z[x]$ with $\deg f(x)\geq 1$. Assume
that $p_1,p_2,\ldots,p_m$ with $p_1<p_2<\ldots<p_m$ are the only
primes of $\Z$ which divide $f(k)$ for any $k\in\Z$ such that
$f(k)\not=0$. Let $f(x)=a_nx^n+\ldots + a_1x+a_0$ and suppose
$a_n>0$. Clearly, $a_0\not=0$. Then we can pick $l$ large enough
so that ${p_i}^l\nmid a_0 = f(0)$ for $i=1,\ldots,m$. Since
$a_n>0$, we can choose $k>l$ such that
${p_m}^{ml+1}<f({p_1}^k{p_2}^k\cdots {p_m}^k)$, but
${p_1}^k{p_2}^k\cdots {p_m}^k$ is an integer, hence by hypothesis
\begin{equation}\label{equation1}
f({p_1}^k{p_2}^k\cdots {p_m}^k) = {p_1}^{j_1}{p_2}^{j_2}\cdots
{p_m}^{j_m},
\end{equation}
for some $j_1,j_2,\ldots,j_m \in \Z^+\cup\{0\}$.\\
Note that ${p_1}^{j_1}{p_2}^{j_2}\cdots {p_m}^{j_m} \leq
{p_m}^{j_1+\ldots+j_m}$, so $f({p_1}^k{p_2}^k\cdots {p_m}^k)\leq
{p_m}^{j_1+\ldots+j_m}$. Hence,
${p_m}^{ml+1}<{p_m}^{j_1+j_2+\ldots +j_m}$. Therefore $ml+1 <
j_1+j_2+\ldots + j_m$ and so for some $i$, $l\leq j_i$. By
$(\ref{equation1})$, we obtain ${p_i}^{l}|f({p_1}^k{p_2}^k\cdots
{p_m}^k) = a_n({p_1}^k{p_2}^k\cdots {p_m}^k)^n+\ldots +
a_1({p_1}^k{p_2}^k\cdots {p_m}^k)+a_0$, therefore ${p_i}^{l}|a_0$,
which is a contradiction.
\end{proof}
The following Corollary provides many \textit{principal ideal domains} ($PID$) that satisfies $IPP$.
\begin{cor}\label{ZnsatisfiesIPP}
For each $n\geq 1$, the ring $\Z\left[\frac{1}{n}\right]$
satisfies $IPP$.
\end{cor}
\begin{proof}
Let $D=\Z\left[\frac{1}{n}\right]$. Let $g(x)\in D[x]$ with $\deg
g(x) \geq 1$. There exists $m\in \Z$ such that $mg(x)\in \Z[x]$. By Lemma \ref{ZsatisfiesIPP}
\[\{p\in P: (\exists k\in \Z)(mg(k)\not=0 \text{ and } p|mg(k)\}\]
is infinite, where $P$ is the set of primes of $\Z$. Therefore
\[\{p\in P: (\exists k\in \Z)(g(k)\not=0 \text{ and } p|g(k)\}\]
is infinite. Hence, if $H=P-\{p\in P: p|n\}$ is the set of primes
of $D$, we obtain that
\[\{p\in H: (\exists k\in D)(g(k)\not=0 \text{ and } p|g(k)\}\]
is infinite. Therefore $D$ satisfies $IPP$.
\end{proof}
The following result generalizes Corollary \ref{ZnsatisfiesIPP}.
\begin{prop}\label{dSDring}
Let $D$ be a $UFD$ and $K= Q(D)$ the quotient field of $D$.
Suppose $D\subseteq S\subseteq K$, where $S$ is a domain, and
suppose $dS\subseteq D$ for some nonzero element $d\in D$. Then
$D$ satisfies $IPP$ if and only if $S$ satisfies $IPP$.
\end{prop}
\begin{proof}
$(\Rightarrow)$. Suppose that $D$ satisfies $IPP$. Note that $S\subseteq
D\left[\frac{1}{d}\right]$. Let $g(x)\in S[x]$ with $\deg g(x) \geq
1$. Because $D$ is a $UFD$, there exists $m\in D$ with $m\not=0$
such that $mg(x)\in D[x]$. Moreover, since $D$ satisfies $IPP$ the set
\[\{p\in P: (\exists k\in D)(mg(k)\not=0 \text{ and } p|mg(k)\}\]
is infinite, where $P$ is the set of primes of $D$. Therefore
\[\{p\in P: (\exists k\in D)(g(k)\not=0 \text{ and } p|g(k)\}\]
is infinite. Note that if $p$ is a prime such that $p|d$ then $p$
is a unit of $D\left[\frac{1}{d}\right]$. Thus, the primes of
$D\left[\frac{1}{d}\right]$ are the primes $p$ in $D$ such that
$p\nmid d$. It follows that the primes in $S$ are the primes $p\in P$
such that $p\nmid d$. Hence, if $P-\{p\in P: p|n\}\supseteq H$,
where $H$ is the set of primes of $S$, we obtain that
\[\{p\in H: (\exists k\in S)(g(k)\not=0 \text{ and } p|g(k)\}\]
is infinite. Therefore $S$ satisfies $IPP$.\\
$(\Leftarrow)$. Suppose that $S$ satisfies $IPP$. Let $f(x)\in D[x]$ with $\deg f \geq 1$.
Assume that $p_1,\ldots,p_m$ are the only primes of
$D$ which divide $f(k)$, for any $k\in D$ such that $f(k)\not=0$.
Define $g(x)=f(dx)$. Note that $g(x)\in S[x]$ and $\deg g(x)\geq 1$.
Let $k\in S$ such that $g(k)\not=0$. Then $g(k) = f(dk)\not=0$.
Also $dS\subseteq D$, so $dk\in D$. Let $p$ be a prime in
$S$ such that $p|g(k)$, then $p=p_i$ for some $i=1,\ldots,m$ because primes in $S$ are also primes
in $D$. This is a contradiction. Therefore $S$ does not satisfy $IPP$.
\end{proof}
\begin{defi}\label{DPP}
A domain $D$ satisfies the \textit{\textbf{degree polynomial property
(DPP)}} if given $g(x),f(x) \in D[x]$ such that for all $k\in D$,
$(g(k)\not=0 \Rightarrow g(k)\vert f(k))$ implies $f(x) = 0$ or
$\deg f(x)\geq \deg g(x)$.
\end{defi}
There is no field $K$ satisfying $DPP$. To see this, take $f(x)=1$ and $g(x)=x$ in $K[x]$. Notice that for all $k\in K$ such that $g(k)\not=0$ we have that $g(k)|f(k)$, however $f(x)\not=0$ and $\deg f(x) < \deg g(x)$.
\begin{exam}
In Section 6, we shall prove that the ring $\Z[W]$, where
\[W:=\{1/p: p \text{ is prime and } p\equiv 1\mod 4 \text { or } p=2 \},\]
does not satisfy $DPP$. The units in this ring are elements $\frac{c}{d}$ with $c\equiv 0 \mod p$ and
$p\equiv 1\mod 4$. It follows that the ring $\Z[W]$ is not a field.
\end{exam}
\begin{lem}\label{ZsatisfiesDPP}
Let $g(x),f(x)\in \Z[x]$ such that $(g(k)\not=0\Rightarrow g(k)|f(k))$, for $k\in \Z$ arbitrary large. Then $f(x)=0$ or $\deg f(x)\geq \deg g(x)$.
\end{lem}
\begin{proof}
Let $g(x)= a_nx^n+\ldots + a_1x+a_0$ and $f(x)=b_mx^m+\ldots + b_1x+b_0$ be polynomials in $\Z[x]$. Without loss of generality, suppose $a_n,b_m > 0$. Assume $(g(k)\not=0 \Rightarrow g(k)|f(k))$, for $k\in \Z$ arbitrary large. If $\deg f(x)=m < n= \deg g(x)$ then (by elementary calculations) we can find $k\in \Z$ large enough such that $g(k)\not=0$ and $a_nk^n+\ldots+a_1k+ a_0 > b_mk^m+\ldots+b_1k+ b_0$. This is a contradiction.
\end{proof}
The following result is an immediate consequence of Lemma \ref{ZsatisfiesDPP}.
\begin{cor}
The ring $\Z$ satisfies $DPP$.
\end{cor}
\begin{prop}\label{D[x]satisfiesDPP}
Let $D$ be a domain. Given $g(y),f(y)\in D[x][y]$ such that for arbitrary large $t$, $g(x^t)|f(x^t)$. Then $f(y) = 0$ or $\deg_y f(y)\geq \deg_y g(y)$.
\end{prop}
\begin{proof}
Let $g(y),f(y)\in D[x][y]$ and suppose $g(x^t)|f(x^t)$ for $t$ arbitrary large. By $\deg_yf(y)$ we mean the highest exponent of $y$ in $f(y)$. Assume that $f(y)\not=0$ and $m=\deg_yf(y) < \deg _yg(y) = n$. Let $g(y) = a_n(x)y^n+\ldots +a_1(x)y+a_0(x)$ and $f(y) = b_m(x)y^m+\ldots +b_1(x)y+b_0(x)$. By hypothesis, $g(x^t)|f(x^t)$, for $t$ arbitrary large, therefore if $h(x)= g(x^t)= a_n(x)x^{tn}+\ldots+a_1(x)x^t+a_0(x)$ and $l(x)=f(x^t)= b_m(x)x^{tm}+\ldots+b_1(x)x^t+b_0(x)$ we have $h(x)|l(x)$. Pick $t$ large enough such that $\deg h(x) = \deg (a_n(x)+ tn)$ and $\deg l(x) = \deg(b_m(x)+tm)$, $f(x^t)\not=0$ and $t>\frac{\deg b_m(x)-\deg a_n(x)}{n-m}$, so $\deg h(x) > \deg l(x)$. Since $h(x)|l(x)$, we obtain $l(x)=0$ or $\deg l(x) \geq \deg h(x)$. In any case we have a contradiction. Therefore $f(y)=0$ or $\deg_y f(y)\geq \deg_y g(y)$.
\end{proof}
The next Corollary shows that a ring of polynomials over any
domain always satisfies $DPP$. Its proof follows from Proposition
\ref{D[x]satisfiesDPP}.
\begin{cor}\label{dpppoly}
Let $D$ be an integral domain. The ring of polynomials $D[x]$ satisfies $DPP$.
\end{cor}
In particular, $\Z[x]$ satisfies $DPP$ and using that $D[x][y] =
D[x,y]$ we have that $\Z[x_1,...,x_n]$ also satisfies $DPP$. Notice that Corollary \ref{dpppoly} also implies that $K[x_1,\ldots,x_n]$ satisfies $DPP$ as well, for any field $K$.
\begin{defi}\label{EPP}
Let $D$ be a \textit{UFD}. $D$ satisfies
the \textit{\textbf{evaluation polynomial property (EPP)}} if given
$f(x),g(x) \in D[x]$ with $g(x)$ primitive, $\deg g(x)\geq 1$ and
for all $k \in D$, $(g(k)\not=0 \Rightarrow g(k)\vert f(k))$, then
$g(x)\vert f(x)$ in $D[x]$. Of course, this is only true when $D$
is infinite (otherwise $D$ is a field).
\end{defi}
There is no an infinite field $K$ satisfying $EPP$. To show this, take $f(x)=1$ and $g(x)=x$ in $K[x]$ where $K$ is an arbitrary infinite field (e.g. $\R$). For all $k\in K$ such
that $g(k)\not=0$ we have that $g(k)|f(k)$ but $g(x)\nmid f(x)$.
On the other hand notice that $5|k^5-k$ for any $k\in\Z$, but certainly
$5\nmid x^5-x$ in $\Z[x]$. This does not prove that the ring of integers does not satisfy $EPP$ (actually it does as we show later), since the constant polynomial $g(x)=5$ is not primitive.
The following Proposition provides a characterization for $EPP$
property.
\begin{prop}\label{EPPirre}
Let $D$ be a \textit{UFD}. $D$ satisfies $EPP$ if and only if given
$f(x)$, $g(x)$ polynomials in $D[x]$ with $g(x)$ irreducible,
$\deg g(x)\geq 1$ and for all $k \in D$, $(g(k)\not=0 \Rightarrow
g(k)\vert f(k))$, then $g(x)\vert f(x)$ in $D[x]$.
\end{prop}
\begin{proof}
See \cite[pg~30]{caceres}.
\end{proof}
Now, we show that in a $UFD$, satisfying $DPP$ is the same
as satisfying $EPP$.

\begin{prop}\label{prop3}
Let $D$ be a $UFD$. $D$ satisfies $DPP$ if and only if $D$ satisfies
$EPP$
\end{prop}
\begin{proof}
$(\Rightarrow)$. Let $D$ be a $UFD$ satisfying $DPP$. Let $f(x),g(x)\in D[x]$
with $g(x)$ primitive, $\deg g(x)\geq 1$ and such that for all
$k\in D$, $(g(k)\not=0 \Rightarrow g(k)|f(k))$. Since $D$
satisfies $DPP$, we obtain $f(x)=0$ or $\deg f(x)\geq \deg g(x)$. If
$f(x)=0$, we are done. Put $g(x)=a_nx^n+\ldots +a_1x+a_0$. By the
usual Division Algorithm, we can find $s\in \Z$ and $q(x),r(x)\in D[x]$
such that
\begin{equation}\label{equation2}
{a_n^s}f(x)=g(x)q(x)+r(x)
\end{equation}
with $\deg r(x) < \deg g(x)$. Since for all $k\in D$, $(g(k)\not=0
\Rightarrow g(k)|f(k))$. Then for all $k\in D$, $(g(k)\not=0
\Rightarrow g(k)|r(k))$. But $D$ satisfies $DPP$, so $r(x)=0$ or
$\deg r(x)\geq \deg g(x)$; thus $r(x)=0$. It follows from
$(\ref{equation2})$ that $g(x)|a_n^sf(x)$. Since $g(x)$ is primitive
and $\deg g(x)\geq 1$, by Gauss' Lemma we obtain
$g(x)|f(x)$. Therefore $D$ satisfies $EPP$.\\
$(\Leftarrow)$. Suppose $D$ satisfies $EPP$. Let $f(x),g(x)\in D[x]$
such that for all $k\in D$, $(g(k)\not=0 \Rightarrow g(k)|f(k))$.
If $\deg g(x)\leq 0$, the result is clear. Suppose $\deg g(x)\geq
1$. Then $g(x)= C(g(x))h(x)$ where $C(g(x))$ is the content of
$g(x)$ and $h(x)$ is a primitive polynomial in $D[x]$ with $\deg
h(x)=\deg g(x)$. By hypothesis, for all $k\in D$, $(h(k)\not=0
\Rightarrow h(k)|f(k))$. Since $D$ satisfies $EPP$ we have $h(x)|f(x)$.
Then $f(x)=0$ or $\deg f(x)\geq \deg h(x)=\deg g(x)$. Therefore
$D$ satisfies $DPP$.
\end{proof}
We obtain the following immediate results from Proposition \ref{prop3} and  Corollary \ref{ZsatisfiesDPP}.
\begin{cor}
The ring $\Z$ satisfies $EPP$.
\end{cor}
\begin{cor}\label{DxsatisfiesEPP}
Let $D$ be a $UFD$. $D[x]$ satisfies $EPP$.
\end{cor}
By Corollary \ref{DxsatisfiesEPP}, we have in  particular that
$\Z[x_1,...,x_n]$ and $K[x_1,\ldots,x_n]$ satisfy $EPP$, where $K$ is any infinite field.
\begin{defi}\label{SEPP}
Let $D$ be a $UFD$. $D$ satisfies the \textit{\textbf{strong
evaluation polynomial property}} $(SEPP)$ if for each $f(x), g(x)
\in D[x]$ where $g(x)$ is irreducible with $\deg g\geq 1$ there
exists $I_{g(x)} \subseteq D$ infinite, such that if $H$ is
infinite and $H\subseteq I_{g(x)}$, then for all $k\in H$,
$(g(k)\not=0 \Rightarrow g(k)\vert f(k))$, implies $g(x)\vert
f(x)$.
\end{defi}
\begin{prop}\label{ZsatisfiesSEPP}
Suppose $f(x),g(x)\in \Z[x]$ with $g(x)$ primitive, $\deg g(x)\geq
1$ and such that $(g(k)\not=0 \Rightarrow g(k)|f(k))$, for $k\in
\Z$ arbitrary large, then $g(x)|f(x)$ in $\Z[x]$.
\end{prop}
\begin{proof}
Let $f(x), g(x)\in Z[x]$ with $g(x)$ primitive, $\deg g(x)\geq 1$
and such that $(g(k)\not=0 \Rightarrow g(k)|f(k))$, for $k\in \Z$
arbitrary large. By Lemma \ref{ZsatisfiesDPP} we obtain that
$f(x)=0$ or $\deg f(x)\geq \deg g(x)$. If $f(x)=0$, we are done.
Suppose $\deg f(x)\geq \deg g(x)$ and let
$g(x)=a_nx^n+a_{n-1}x^{n-1}+\ldots+a_0$. By the usual Division
Algorithm, we can find $s\in \Z$ and $q(x),r(x)\in \Z[x]$ such that
$a_n^sf(x)=g(x)q(x)+r(x)$ with $\deg r(x)<\deg g(x)$. Since
$(g(k)\not=0 \Rightarrow g(k)|f(k))$ for $k$ arbitrary large, then
\[(g(k)\not=0 \Rightarrow g(k)|r(k)),\]
for $k$ arbitrary large. By Lemma \ref{ZsatisfiesDPP}, $r(x)=0$ or
$\deg r(x)\geq \deg g(x)$. Therefore $r(x)=0$, which implies that
$g(x)|a_n^sf(x)$ with $g(x)$ primitive and $\deg g(x)\geq 1$.
By Gauss' Lemma,  $g(x)|f(x)$ in $\Z[x]$.
\end{proof}
\begin{cor}
$\Z$ satisfies $SEPP$.
\end{cor}
\begin{proof}
Let $g(x)\in \Z[x]$, irreducible with $\deg g(x)\geq 1$. Let
$I_{g(x)}=\Z^{+}$. The result now follows from Proposition \ref{ZsatisfiesSEPP}.
\end{proof}
The following Proposition provides examples of domains satisfying
$EPP$.
\begin{prop}\label{prop4}
Let $D$ be a domain. If $D$ satisfies $SEPP$, then $D$ satisfies
$EPP$.
\end{prop}
\begin{proof}
Suppose $D$ satisfies $SEPP$. Let $f(x),g(x)\in D[x]$, with $g(x)$
primitive and $\deg g(x)\geq 1$. Suppose that
\begin{equation}\label{equation3}
\text{for all $k\in D$, $(g(k)\not=0 \Rightarrow g(k)|f(k))$.}
\end{equation}
Actually, by Proposition \ref{EPPirre}, we can assume that $g(x)$
is irreducible. By hypothesis, there exists $I_{g(x)}\subseteq D$
infinite, such that
\begin{align}\label{equation4}
\text{for each $H\subseteq I_{g(x)}$ infinite,}\\
\text{if for each $k\in H$, $(g(k)\not=0 \Rightarrow g(k)|f(k))$,
then $g(x)|f(x)$.}\\ \notag
\end{align}
By $(\ref{equation3})$ we have that for all $k \in I_{g(x)}$,
$(g(k)\not=0 \Rightarrow g(k)|f(k))$. In particular, for
$H=I_{g(x)}$ in $(\ref{equation4})$, we obtain $g(x)|f(x)$.
Therefore $D$ satisfies $EPP$.
\end{proof}
The following Proposition says that in a $UFD$, $IPP$ implies $SEPP$. Its proof uses
ultraproducts, which is a topic not related to the theory of this
paper.
\begin{prop}\label{prop5}
Let $D$ a $UFD$. If $D$ satisfies $IPP$ then $D$ satisfies $SEPP$.
\end{prop}
\begin{proof}
See \cite[pg~36]{caceres}.
\end{proof}
\begin{prop}\label{prop10}
Let $D$ be a $UFD$ with at least one prime and with finitely many
units, then $D$ satisfies $EPP$.
\end{prop}
\begin{proof}
See \cite[pg~38]{caceres}
\end{proof}
The converse of Proposition \ref{prop10} is not true in general.
The ring $\Z\left[\frac{1}{n}\right]$ satisfies $EPP$ by
Proposition \ref{ZnsatisfiesIPP}, Proposition \ref{prop5} and
Proposition \ref{prop4}, but it has infinitely many units; in fact
the units of $\Z\left[\frac{1}{n}\right]$ are the integers $p^j$
with $p$ prime and such that $p|n$. However this ring also
satisfies $DPP$ and $SEPP$.
\section{$D$-rings}
\begin{defi}\label{D-Ring}
Let $D$ be a domain and $K=Q(D)$ its quotient field. $D$ is a \textit{\textbf{$D$-ring}} if given
$f(x),g(x) \in D$ such that, if for almost all $k\in D$, $g(k)\vert f(k)$,
then $\frac{f(x)}{g(x)}\in K[x]$
\end{defi}
A field is never a $D$-ring. To see this, let $K$ be a field. Take $f(x)=x$ and $g(x)=1$, for almost all $k\in D$ we have $f(k)|g(k)$ in $K$ but $\frac{g(x)}{f(x)}\not\in Q(K)[x]=K[x]$.
As we show later, the $D$-ring property is related with rational functions $r(x)$ over $D$ and
polynomials $p(x)$ over $K$ where $K$ its the quotient field of
$D$, such that $r(D),p(D)\subseteq D$. Many interesting results follows from
the $D$-ring property (see \cite[pgs~61-66]{narki} and \cite{hiroshi}).
Our main goal in this section is to show that the $D$-ring property is equivalent to some of the
divisibility properties studied in the previous section.
\begin{lem}\label{ZisaDring}
Let $f(x)$ and $g(x)\in \Z[x]$ such that, for almost all $k\in \Z$, $g(k)|f(k)$. Then $\frac{f(x)}{g(x)}\in \Q[x]$.
\end{lem}
\begin{proof}
If $g(x)$ is a constant-nonzero polynomial, we are done. Assume $\deg g(x)\geq 1$. Let $A=\{k_1,\ldots,k_n\}$ such that for all $k\in \Z-A$, $g(k)|f(k)$. Let $k_1,\ldots,k_s \in A$ such that $g(k_i)\not= 0$ for $i=1,\ldots,s$ and let $\beta = g(k_1)\cdots g(k_s)$. If $s=0$, let $\beta =1$. Then for all $k\in \Z$ such that $g(k)\not=0$, $g(k)|\beta f(k)$. Since $\Z$ satisfies $EPP$ we have that $g(x)|\beta f(x)$ in $\Z[x]$. Hence, there exists $p(x)\in \Z[x]$ such that $\beta f(x)=p(x)g(x)$. So $\frac{f(x)}{g(x)}= \beta^{-1}p(x)\in \Q[x]$.
\end{proof}
We have the following Corollary of Lemma \ref{ZisaDring}.
\begin{cor}\label{ZisDring}
$\Z$ is a $D$-ring.
\end{cor}
Note that by Corollary \ref{ZisDring}, given $f(x)$ and $g(x)$
polynomials with coefficients in $\Z$ such that $g(k)|f(k)$ for almost all $k\in \Z$, implies the existence of a polynomial $h(x)=
\frac{f(x)}{g(x)}\in \Q[x]$ with $h\left(\Z\right)\subseteq
\Z$. For example, if $p$ is a prime in $\Z$, we have that for
any $k\in \Z$, $p|k^p-k$ which implies $\frac{x^p-x}{p}\in\Q[x]$.
\begin{exam}
In the Section 6, we show that the ring $\Z[W]$, where
\[W:=\{1/p: p \text{ is prime and } p\equiv 1\mod 4 \text { or } p=2 \},\]
is not a $D$-ring. We have already shown that this ring is not a field.
\end{exam}
\begin{defi}
Let $D$ be a domain. For any polynomial $f(x) \in D[x]$ denote $S(f)$
the set of all non-zero prime ideals $\mathfrak{P}$ of $D$ such
that the congruence $f(x)\equiv 0\mod\mathfrak{P}$ is solvable in
$D$. This is: there exists $k\in D$ such that $f(k)\in \Prime$. In
particular, if $c\in D$, $S(c)$ is precisely the set of prime
ideals of $D$ that contain $c$.
\end{defi}
\begin{prop}\label{dringeq}
Let $D$ be a domain, $K$ the quotient field of $D$ and $D^\times$
the set of units of $D$ . The following properties are equivalent:
\begin{enumerate}
\item[(1)] $D$ is a $D-ring$.
\item[(2)] Every polynomial over $D$ which satisfies $f(k)\in D^\times$ for almost all $k\in D$ must be a constant.
\item[(3)] For any non-constant polynomial $f(x)\in D[x]$, the set $S(f)$ is non-empty.
\item[(4)] For any non-constant polynomial $f(x) \in D[x]$ and any non-zero $c\in D$, the set $S(f)-S(c)$ is infinite.
\end{enumerate}
\end{prop}
\begin{proof}
See \cite[pgs~61-62]{narki} or \cite[pgs~290-291]{hiroshi}.
\end{proof}
Proposition \ref{dringeq} gives us a very useful tool for proving
results about $D$-rings. The following Corollary gives a characterization of
the $D$-ring property for domains that are not fields, its proof is an immediate consequence of Proposition \ref{dringeq}.
\begin{cor}\label{funit}
Let $D$ be a ring that is not a field and $D^\times$ be the set of
units of $D$. $D$ is not a $D$-ring if and only if there exists a nonconstant
polynomial $f(x)\in D[x]$ such that $f(D)\subseteq D^\times$.
\end{cor}
The following result gives a relation between a $D$-ring and its
\textit{Jacobson Radical} (denoted by $\J(D)$ for any ring $D$).
\begin{prop}\label{jacobsondring}
Let $D$ be a ring that is not a field. If $\mathfrak{J}(D)\not=(0)$ then $D$ is not a $D$-ring.
\end{prop}
\begin{proof}
If $\mathfrak{J}(D)\not=(0)$, then let $c\in \mathfrak{J}(D)$ with
$c\not=0$. We have that the polynomial $f(x)= 1-cx$ satisfies
$f(D)\subseteq D^\times$. By Corollary \ref{funit}, $D$ is not
a $D$-ring.
\end{proof}
There is a relation between $IPP$ and the $D$-ring property. The
$IPP$ talks about infinitely many prime elements, while the
$D$-ring property talks about infinitely many prime ideals. So, in
a $PID$ it is trivial that $IPP$ and the $D$-ring property are
equivalent properties.
Now, we show that any $UFD$ satisfying the $D$-ring property, also satisfies $IPP$.
\begin{prop}\label{prop1}
Let $D$ be a $UFD$. If $D$ is a $D$-ring, then $D$ satisfies $IPP$.
\end{prop}
\begin{proof}
Let $g(x)\in D[x]$ with $\deg g(x)\geq 1$. Suppose that there exists $a\in D$ with
$g(a)=0$. Then, there exists $m\in D$ and $h(x)\in D[x]$ such that
$mg(x)=(x-a)h(x)$. Let $p$ be a prime of $D$ such that $p\nmid m$
and $h(p+a)\not=0$. Note that $D$ has infinitely many primes satisfying this condition. Therefore $mg(p+a)=ph(p+a)$, so $p\vert
mg(p+a)$. By our choice of $p$, we have that $p\vert g(p+a)$.
Therefore the set
\[\{p\in P: (\exists k\in D)(g(k)\not=0 \text{ and } p\vert g(k)\},\]
where $P$ is the set of primes of $D$ is infinite. So, $D$ satisfies $IPP$.
Suppose that $g(a)\not=0$ for all $a\in D$. Assume that $p_1,\ldots,p_n$ are the only primes of $D$ which divide
$g(k)$ for any $k\in D$ such that $g(k)\not=0$. Let $m=p_1\cdots
p_n$. Since $D$ is a $D$-ring the set $S(g)-S(m)$ is not empty. Let
$\mathfrak{P}\in S(g)-S(m)$, then there exists
$k_{\mathfrak{P}}\in D$ such that
$g(k_{\mathfrak{P}})\in\mathfrak{P}$ and $m\not\in\mathfrak{P}$.
By our assumption
\[g(k_\mathfrak{P})= u{p_1}^{m_{1}}{p_2}^{m_{2}}\cdots{p_n}^{m_{n}},\]
where $u\in D^\times$ and $m_i$ is a non-negative integer for $i=1,\ldots,n$. Since $g(k_\mathfrak{P})\in \mathfrak{P}$, then $u\in \mathfrak{P}$ or there exists $j \in\{1,\ldots,n\}$ such that ${p_j}^{m_j}\in\mathfrak{P}$. If $u\in \mathfrak{P}$ then $\mathfrak{P}=D$ and this contradicts that $\mathfrak{P}$ is a prime ideal of $D$. If ${p_j}^{m_k}\in\mathfrak{P}$, then $p_j\in \mathfrak{P}$, therefore $m\in \mathfrak{P}$, and this is also a contradiction.
Therefore $D$ satisfies $IPP$.
\end{proof}
The converse of the previous result is also true, but we need some previous results in order to prove it. The following Proposition shows that domains that satisfies $DPP$
are $D$-rings and viceversa.
\begin{prop}\label{DringiffDPP}
Let $D$ be a domain. $D$ is a $D$-ring if and only if $D$ satisfies $DPP$.
\end{prop}
\begin{proof}
$(\Rightarrow)$. Let $g(x),f(x) \in D[x]$ such that for all $k\in D$, $(g(k)\not=0\Rightarrow g(k)\vert f(k))$. So, $g(k)\vert f(k)$ for almost $k \in D$. Since $D$ is a $D$-ring, then $\frac{f(x)}{g(x)}\in K[x]$, where $K$ is the quotient field of $D$. Therefore, there exists $p(x)\in K[x]$ such that $f(x)=p(x)g(x)$. Suppose that $f(x)\not= 0$, so $\deg f(x) = \deg (p(x)g(x)) = \deg p(x) + \deg g(x)\geq \deg g(x)$, then $D$ satisfies $DPP$.\\
$(\Leftarrow)$. Let $g(x),f(x)\in D[x]$ such that for almost all $k\in D$, $g(k)\vert f(k)$. Let $A=\{k_1,\ldots,k_n\}$ be a finite subset of $D$ such that $g(k)\vert f(k)$ for all $k \in D-A$. Let $k_1,\ldots,k_s \in A$ such that $g(k_i)\not= 0$ for $i=1,\ldots,s$ and let $\beta = g(k_1)\cdots g(k_s)$. If $s=0$, let $\beta = 1$. Then, for all $k\in D$ such that $g(k)\not= 0$ we obtain that  $g(k)\vert \beta f(k)$. Since $D$ satisfies $DPP$, then $\beta f(x)=0$ or $\deg \beta f(x)\geq \deg g(x)$. If $\beta f(x)=0$, then $f(x)=0$, so $\frac{f(x)}{g(x)}\in K[x]$. Suppose that $\deg \beta f(x)\geq \deg g(x)$ and assume $g(x)= a_nx^n+\ldots+a_0$. By The Division Algorithm there exist $q(x),r(x) \in K[x]$ and $s\in D$ such that
\[a_n^s\beta f(x)= g(x)q(x)+r(x),\]
with $r(x)=0$ or $\deg r(x) < \deg g(x)$ and let $\alpha=a_n^s\beta$. Suppose that $\deg r(x)<\deg g(x)$. Then for all $k\in D$ such that $g(k)\not = 0$ implies that $g(k)\vert\alpha f(k)$ and $g(k)\vert g(k)q(k)$. So $g(k)\vert r(k)$. Hence, using again that $D$ satisfies $DPP$ we obtain $r(x) = 0$ or $\deg r \geq \deg g$. Hence $r(x) = 0$ and we obtain that $\alpha f(x) = g(x)q(x)$.
Therefore $\frac{f(x)}{g(x)} = \alpha^{-1}q(x)\in K[x]$. In others words, $D$ is a $D$-ring.
\end{proof}
The following Proposition shows that $UFD$'s satisfying $EPP$
are $D$-rings and viceversa.
\begin{prop}\label{prop2}
Let $D$ be a $UFD$. $D$ is a $D$-ring if and only if $D$ satisfies $EPP$.
\end{prop}
\begin{proof}
$(\Rightarrow)$. Let $f(x),g(x) \in D[x]$ with $g$ primitive and $\deg g(x)\geq 1$ such that for all $k \in D$, $g(k)\not=0 \Rightarrow g(k)\vert f(k)$. It is clear that for almost all $k\in D$, $g(k)\vert f(k)$. Since $D$ is a $D-ring$ we have that
\[\frac{f(x)}{g(x)}=p(x)\in K[x],\]
where $K = Q(D)$ is the quotient field of $D$. Let
\[p(x) = \frac{r_n}{s_n}x^n+\frac{r_{n-1}}{s_{n-1}}x^{n-1}+\ldots+\frac{r_1}{s_1}x+\frac{r_0}{s_0},\]
where $r_i,s_i \in D$, with $s_i \not= 0$ for all $i = 0,\ldots,n$.
Let $m = \lcm(s_n,\ldots,s_0)$ (this element exists, because $D$ is a $UFD$), therefore  $mp(x)\in D[x]$. Take $h(x)=mp(x)$.
Now, we have that
\[mf(x) = mp(x)g(x)=h(x)g(x),\]
with $g(x)$ primitive. By Gauss' Lemma, there exists $q(x)\in
D[x]$ such that $h(x) = mq(x)$, and so
\[mf(x) = mq(x)g(x).\]
Therefore $f(x) = q(x)g(x)$, with $q(x) \in D[x]$; i.e. $g(x)\vert f(x)$ in $D[x]$. Hence, $D$ satisfies $EPP$.\\
$(\Leftarrow)$.  Let $f(x),g(x) \in D[x]$ such that for almost all $k\in D$ we have that $g(k)\vert f(k)$. Let $A=\{k_1,\ldots,k_n\}$ be a finite subset of $D$ such that $g(k)\vert f(k)$ for all $k \in D-A$. Let $k_1,\ldots,k_s \in A$ such that $g(k_i)\not= 0$ for $i=1,\ldots,s$ and let $\beta = g(k_1)\cdots g(k_s)$. If $s=0$, let $\beta = 1$. Then for all $k\in D$ such that $g(k)\not= 0$ we have $g(k)\vert \beta f(k)$. Let $K=Q(D)$ be the quotient field of $D$. We can write $g(x)=\alpha h(x)$ where $h(x)$ is primitive with $\deg h=\deg g \geq 1$ and $\alpha$ is the content of $g(x)$.
Let $k \in D$ such that $h(k)\not=0$. Therefore $g(k)\not=0$ and $g(k)\vert \beta f(k)$; but $h(k)\vert g(k)$, so $h(k)\vert\beta f(k)$. Since $D$ satisfies $EPP$, we have that $h(x)\vert\beta f(x)$ in $D[x]$. Hence, there exists $p(x)\in D[x]$ such that $\beta f(x)=p(x)h(x)$ and so
\[\alpha\beta f(x)= p(x)(\alpha h(x)) = p(x)g(x).\]
Therefore $f(x) = {(\alpha\beta)}^{-1}p(x)g(x)$ where ${(\alpha\beta)}^{-1}p(x) \in K[x]$, i.e. $g(x)\vert f(x)$ in $K[x]$. Hence, $D$ is a $D$-ring.
\end{proof}
\begin{cor}\label{DxisDring}
Let $D$ be a domain. The ring $D[x]$ is a $D$-ring.
\end{cor}
\begin{proof}
Immediate from Proposition \ref{prop2} and Corollary
\ref{DxsatisfiesEPP}.
\end{proof}
Using Corollary \ref{DxisDring} we have that the rings $\Z[x_1,\ldots,x_n]$ and $K[x_1,...,x_n]$,
where $K$ is a field are $D$-rings. Note that by Corollary
\ref{DxisDring} and Proposition \ref{jacobsondring}, we obtain that for any domain
$D$, $\J(D[x])=\{0\}$, for instance,
$\J(\Z[x_1,\ldots,x_n])=\{0\}$.
The ring $\Z$ satisfies all our divisibility properties as well as the ring $D[x_1,\ldots,x_n]$ for any domain $D$. The
following Theorem says that in any $UFD$, the properties $IPP$, $DPP$, $EPP$, $SEPP$ and the $D$-ring property are equivalent.
\begin{thm}\label{equivalence}
Let $D$ be a $UFD$. The following properties are equivalent:
\begin{itemize}
\item[(1)] $D$ is a $D$-ring.
\item[(2)] $D$ satisfies $IPP$.
\item[(3)] $D$ satisfies $DPP$.
\item[(4)] $D$ satisfies $EPP$.
\item[(5)] $D$ satisfies $SEPP$.
\end{itemize}
\end{thm}
\begin{proof}
$(1)\Rightarrow (2)$ from Proposition \ref{prop1}, $(2)\Rightarrow
(5)$ from Proposition \ref{prop5}, $(5)\Rightarrow(4)$ from
Proposition \ref{prop4}, $(4)\Rightarrow(1)$ from Proposition
\ref{prop2} and $(3)\Leftrightarrow(4)$ from Proposition
\ref{prop3}.
\end{proof}
The following Corollary gives infinitely many $PID$'s that are $D$-rings. It is a consequence of Proposition \ref{ZnsatisfiesIPP} and Theorem \ref{equivalence}
\begin{cor}
For all $n\geq 1$, $\mathbb{Z}[\frac{1}{n}]$ is a $D-$ring.
\end{cor}
By Theorem \ref{equivalence} and Corollary \ref{DxisDring}
we have that $D[x]$ with $D$ a domain, satisfies all
divisibility properties $IPP$, $DPP$, $EPP$ y $SEPP$. Furthermore
$D[x]$ is also a $D$-ring. Therefore, we obtain a number of
rings satisfying our divisibility properties, for example:
$\Z[x_1,\ldots,x_n]$, $\Z_p[x_1,\ldots,x_n]$ where $p$ is an
integer prime and the ring $\R[x_1,\dots,x_n]$.

\begin{cor}\label{dSDpp}
Let $D$ be a $UFD$ and $K= Q(D)$ be the quotient field of $D$. Suppose
$D\subseteq S\subseteq K$, where $S$ is a domain, and suppose
$dS\subseteq D$ for some nonzero element $d\in D$. Then $D$
is a $D$-ring (resp. satisfies $DPP$, $EPP$ or $SEPP$) if and only if $S$ is a $D$-ring (resp. satisfies $DPP$, $EPP$ or $SEPP$).
\end{cor}
\begin{proof}
Easy from Proposition \ref{dSDring} and Theorem \ref{equivalence}.
\end{proof}
We will assume the following results proven in \cite[pg~299]{hiroshi}.
\begin{prop}
Suppose $D$ is a domain such that $\Z\subseteq D\subseteq \Q$. If $D$ is a non-$D$-ring, then so is every ring between $D$ and $\Q$. If $D$ is a $D$-ring, then so is every ring between $\Z$ and $D$.
\end{prop}
\begin{prop}
Among the subdomains of $\Q$ that are infinitely generated over $\Z$, there are infinitely many $D$-rings and infinitely many non-$D$-rings.
\end{prop}
In the following example it is necessary to know results from
Algebraic Number Theory, topic far away from the theory in this
paper.
However, the reader could find more details in
\cite[pg~293]{hiroshi}.
\begin{exam}
Let $V$ be a set of rational primes $p$ such that $\sum_{p\in V}
1/p$ converges. Let $U$ be the set of all $p^{-1}$ $(p\in V)$.
Then $S=\Z[U]$ is a $D$-ring.
\end{exam}
Note that $\Z[U]$ is a infinitely generated ring over $\Z$
contained in $\Q$.
\section{Infinitely Many Primes}
A result that is interesting is the following:
\begin{prop}\label{infiniteprimes1}
Let $D$ be a $UFD$ with at least one prime and finitely many
units, then $D$ has infinitely many primes.
\end{prop}

\begin{proof}
By Proposition \ref{prop10}, $D$ satisfies $EPP$; therefore $D$
satisfies $IPP$. Then $D$ has infinitely many primes.
\end{proof}
We shall give a direct proof of the previous Proposition but before that we need to prove some Lemmas first.
\begin{lem}[Kaplanski]\label{infmax}
Let $D$ be an infinite domain with a finite number of units, then
$D$ has an infinite number of maximal ideals.
\end{lem}
\begin{proof}
Suppose that $D$ has a finite number of maximal ideals
$\M_1,...,\M_n$. Then the Jacobson Radical of $D$ is $\J(D) =
\bigcap_{k=1}^{n}\M_k$. Because $\M_k\not=(0)$ for all
$k=1,...,n$, then there exists $m_k\in \M_k$ with $m_k \not=0$ for
each $k=1,...,n$. Therefore $m = m_1\cdots m_n \in \M_1\cdots \M_n
\subseteq \J(D)$ with $m\not=0$, hence $\J(D)\not=(0)$. Let
$r\in\J(D)$ with $r\not=0$, then $1-r$ is a unit. Let $U=
\{u_1,...,u_s\}$ the set of units of $D$, then $r=1-u_i$ for some
$i=1,...,s$; therefore $\J(D)$ is finite.
Let $x\in \J(D)$, since $\J(D)$ is finite then for all $n\geq 1$,
there exists $k\leq n$ such that $x^n= x^k$, so $x^{n-k}=1$,
therefore $1\in \J(D)$. Then we have that $\J(D)=D$, so $D$ is
finite, contradicting that $D$ is infinite.
\end{proof}
\begin{lem}\label{primeunions}
Let $\Prime_1,\Prime_2,\ldots,\Prime_n$ be prime ideals of a domain $D$ and let
$\A$ be an ideal of $D$ contained in $\bigcup_{i=1}^n \Prime_i$. Then
$\A\subseteq \Prime_i$ for some $i$ with $i=1,\ldots,n$.
\end{lem}
\begin{proof}
See \cite[pg~8]{atiyah}.
\end{proof}
Now we prove a stronger result than Proposition \ref{infiniteprimes1}. Actually, we could say that the following result is a generalization of Euclid's Theorem about primes.
\begin{prop}\label{infiniteprimes2}
Let $D$ be an infinite $UFD$ with a finite number of units, then $D$ has an infinite number of primes.
\end{prop}
\begin{proof}
Suppose that $p_1,p_2,\ldots,p_n$ are the unique primes in $D$. Let $D^\times$ be the multiplicative group of $D$; $\Gamma=\{\LA p_1 \RA,\ldots,\LA p_n \RA\}$ and $S$ be the set of all maximal ideals of $D$. Since $D$ is a $UFD$ we have that
\[D=\LA p_1 \RA \cup \LA p_2 \RA \cup \dots \cup \LA p_n \RA \cup D^\times.\]
We claim that $S\subseteq\Gamma$.
Let $\M\in S$, then $\M\subseteq D$. Hence $\M\subseteq \LA p_1 \RA\cup \LA p_2 \RA\cup \dots \cup \LA p_n \RA$, where $\LA p_i \RA$ is a prime ideal of $D$ for $i=1,\ldots,n$. By Lemma \ref{primeunions}, we have that $\M\subseteq \LA p_{i_0} \RA$ for some $i_0\in \{1,\ldots,n\}$. But $\M$ is a maximal ideal of $D$, so $\M=\LA p_{i_0} \RA$. Then $\M \in \Gamma$. This proves that $S\subseteq \Gamma$. But $\Gamma$ is a finite set and by Lemma \ref{infmax}, $S$ should be infinite. This is a contradiction. Therefore $D$ has an infinite number of prime elements.
\end{proof}
It is clear that Proposition \ref{infiniteprimes1} is a direct
consequence of Proposition \ref{infiniteprimes2}.
It follows from Proposition \ref{infmax} that if $D$ is an infinite $PID$ with a finite number of units, then $D$ has an infinitely many prime elements.
\section{Many Variables}
The following result shows that we can generalize our divisibility properties of polynomials in one variable to polynomials in two variables. Since we can extend the same argument to polynomials in arbitrary number of variables, it is sufficient to show the two variables case only.
\begin{prop}\label{twovariablesdpp}
Let $D$ be a domain. $D$ satisfies $DPP$ if and only if given $f(x,y), g(x,y)\in D[x,y]$ such that for all $a, b\in D$, $(g(a,b)\not=0 \Rightarrow g(a,b)|f(a,b))$ then $f(x,y)=0$ or $\deg_y f(x,y)\leq \deg_y g(x,y)$. Note that we can replace $\deg_y$ by $\deg_x$.
\end{prop}
\begin{proof}
$(\Leftarrow)$. Since $D[x]\subseteq D[x,y]$ this implication is clear.\\
$(\Rightarrow)$. Suppose that $f(x,y)\not=0$. Let $g(x,y)=c_n(x)y^n+\cdots+c_1(x)y+x_0(x)$ and $f(x,y)=b_m(x)y^m+\cdots+b_1(x)y+b_0(x)$ with $c_n(x), b_m(x)\not=0$. Let $a\in D$ such that $c_n(a), b_m(a)\not=0$, i.e. $f(a,y)\not=0$. Define $h(y)=g(a,y)$ and $l(y)=f(a,y)$. Note that $h(y), l(y)\in D[y]$ and $\deg h(y)=\deg_y g(x,y)$ and $\deg l(y)=\deg_y f(x,y)$. Let $b\in D$ such that $h(b)=g(a,b)\not=0$. By hypothesis, $h(b)=g(a,b)|f(a,b)=l(b)$. Since $D$ satisfies $DPP$, we have that $l(y)=0$ or $\deg h(y)\leq \deg l(y)$. If $l(y)=0$ then $f(a,y)=0$, contradicting that $f(a,y)\not=0$. Then $\deg_y g(x,y)=\deg h(y) \leq \deg l(y)=\deg_y f(x,y)$.
\end{proof}
We have the following Corollary from Proposition \ref{twovariablesdpp}.
\begin{prop}
Let $D$ be a domain. $D$ satisfies $DPP$ if and only if given $f(x_1,\ldots,x_n), g(x_1,\ldots,x_n)\in D[x_1,\ldots,x_n]$ such that for all $a_1\ldots,a_n\in D$,
\[g(a_1,\ldots,a_n)\not=0 \Rightarrow g(a_1,\ldots,a_n)|f(a_1,\ldots,a_n).\]
Then $f(x_1,\ldots,x_n)=0$ or $\deg_{x_i} f(x_1,\ldots,x_n)\leq \deg_{x_i} g(x_1,\ldots,x_n)$ for all $i=1\ldots, n$.
\end{prop}
\begin{cor}
Let $D$ be a $UFD$. $D$ satisfies $EPP$ if and only if given $f(x,y), g(x,y)\in D[x,y]$ with $g(x,y)$ primitive with respect to the variable $y$ and $\deg_y g(x,y)\leq 1$, such that for all $a,b\in D$, $(g(a,b)\not=0\Rightarrow g(a,b)|f(a,b))$ then $g(x,y)|f(x,y)$.
\end{cor}
\begin{proof}
$(\Leftarrow).$ Since $D[x]\subseteq D[x,y]$, this implication is clear.\\
$(\Rightarrow).$ Suppose that $D$ satisfies $EPP$. By Theorem \ref{equivalence}, $D$ also satisfies $DPP$. It follows from Proposition \ref{twovariablesdpp} that $f(x,y)=0$ or $\deg_y g(x,y)\leq \deg_y f(x,y)$. Let $g(x,y)=c_n(x)y^n+\cdots+c_1(x)y+c_0(x)$. By the usual Division Algorithm, we can find $s\in \Z$ and $q(x,y), r(x,y)\in D[x,y]$ such that
\begin{equation}\label{onemoreequation}
c_n^s(x)f(x,y)=q(x,y)g(x,y)+r(x,y),
\end{equation}
with $r(x,y)=0$ or $\deg_y r(x,y) <\deg_y g(x,y)$. Since for all $a,b\in D$ $(g(a,b)\not=0 \Rightarrow g(a,b)|f(a,b))$, then for all $a,b\in D$ $(g(a,b)\not=0\Rightarrow g(a,b)|r(a,b))$. Since $D$ satisfies $DPP$, by Proposition \ref{twovariablesdpp} we obtain $r(x,y)=0$ or $\deg_y g(x,y)\leq \deg_y r(x,y)$. Thus $r(x,y)=0$. By $(\ref{onemoreequation})$, $g(x,y)|c_n^s(x)f(x,y)$. Since $g(x,y)$ is primitive with respect to the variable $y$ and $\deg_y g(x,y) \geq 1$, by Gauss' Lemma we obtain that $g(x,y)|f(x,y)$.
\end{proof}

\section{$Int(D)$}

\begin{defi}
Let $D$ be a domain and $K$ be its quotient field. The set $Int(D)$ is the ring of all polynomials $p(x)$ in $K[x]$, such that $p(D)\subseteq D$.
\end{defi}

We have that $D[x]\subseteq Int(D) \subseteq K[x]$.

For example: for any prime $p$, the polynomial
$f(x)=\frac{x^p}{p}-\frac{x}{p}\in Int(\Z)$ because $f(x)\in
\Q[x]$ and $f(\Z)\subseteq \Z$.

\begin{defi}
Let $D$ be a domain. The set $\S(D)$ is the ring the all rational
functions of $D(x)$ such that, given $r(x)\in \S(D)$, for all
$k\in D$ with $k$ in the domain of $r(x)$ implies that $r(k)\in D$.
\end{defi}

For example, for $n>1$, $r(x)=\frac{1-x^n}{1-x}\in \S(\Z)$. In the
next Section we will give no trivial examples of polynomials
$f(x)$ and $g(x)$ with coefficients in $\Z$ such that for almost
all $k\in \Z$, $g(k)|f(k)$ implies $g(x)|f(x)$.

We always have that $Int(D) \subseteq \S(D)$. But if $K$ is a
field $\S(K)\not\subseteq Int(K)$, because $r(x)=\frac{1}{x}\in
\S(K)$, but $r(x)\not\in Int(K)$.

We give an alternative characterization of the divisibility property $EPP$.

\begin{prop}
Let $D$ be a $UFD$. $D$ satisfies $EPP$ if and only if given $f(x),g(x) \in D[x]$ with
$\deg g \geq 1$ such that $\frac{f(x)}{g(x)} \in \S(D)$ then $g(x)\vert f(x)$ in $D[x]$.
\end{prop}
The following Proposition provides  a characterization of $D$-rings.
\begin{prop}\label{SD=IntD}
Let $D$ be a domain. $D$ is a $D$-ring if and only if $\S(D)=Int(D)$.
\end{prop}
\begin{proof}
See \cite{narki}.
\end{proof}
Note that by Proposition \ref{SD=IntD} and the fact that $\Z$ is
a $D$-ring we have that for any polynomial $h(x)\in \Q[x]$ with $h(\Z)\subseteq \Z$, there exist polynomials $f(x), g(x) \in
\Z[x]$ such that $h(x)=\frac{f(x)}{g(x)}$.
\begin{exam}
There are no localizations $\mathbb{Z}_{(p)}$ of $\Z$ with
respect to a prime $p$ being $D$-rings. In fact, define
$r(x)=\frac{1}{1+px}$. Let $\alpha \in \mathbb{Z}_{(p)}$, then
$\alpha = \frac{a}{b}$ with $a,b \in \mathbb{Z}$ and $b\not\in
(p)$. Then $r(\alpha)= \frac{b}{b+ap}$. It's clear that
$b+ap\not\in (p)$, so $r(\alpha)\in \mathbb{Z}_{(p)}$. Therefore
$r(x)\in \S(\mathbb{Z}_{(p)})$, but $r(x)\not\in
Int(\mathbb{Z}_{(p)})$. Hence $\mathbb{Z}_{(p)}$ is not a
$D$-ring.
\end{exam}

\section{Examples}

In the first part of this section we give nontrivial examples of polynomials with
coefficients in $\Z$ such that for almost all $k\in \Z$
$g(k)|f(k)$ implies that $g(x)|f(x)$ in $\Z[x]$. In the second
part we give a nontrivial ring generated over $\Z$
contained in $\Q$ that is not a $D$-ring.\\
\subsection{Pell's equation}
Consider the following equation:
\begin{equation}\label{pell'sequation}
x^2-dy^2=1,
\end{equation}
where $d$ is a integer that is not a square. Equation $(\ref{pell'sequation})$ is
named as \textit{Pell's equation}. Lagrange proved that
$(\ref{pell'sequation})$ has an infinite number of nontrivial integer solutions (see \cite[pg~320]{burton}).
We are interested on studying a particular case of $(\ref{pell'sequation})$:
\begin{equation}\label{pell}
x^2-(a^2-1)y^2=1,
\end{equation}
where $a\in \mathbb{Z}-\{0,-1\}$.
In \cite{matiya} it is proved the
following recursive formula describing all solution of $(\ref{pell})$.
These are also known as \textit{Lucas' sequences} : if $|a| \geq 2$:
\begin{align}\label{lucas1}
X_0(a)&=1,& X_1(a)&=a,& X_{n+1}(a)&=2aX_n(a)-X_{n-1}(a);\\
Y_0(a)&=0,& Y_1(a)&=1,& Y_{n+1}(a)&=2aY_n(a)-Y_{n-1}(a).
\end{align}
If $a=1$, define for all $n\geq 0$:
\begin{align}\label{lucas2}
X_{n}(1)=1,\\
Y_{n}(1)=n.
\end{align}
Table \ref{table} shows the values for $X_a(n)$ and $Y_a(n)$ with $|a| \geq 2$ for $n=0,1,\ldots,8.$\\
\begin{table}[ht]
\centering
\begin{tabular}{|r|r|r|}\hline
$n$ & $X_{n}(a)$ & $Y_{n}(a)$\\
\hline $0$ & $1$ & $0$ \\
\hline $1$ & $a$ & $1$ \\
\hline $2$ & $2a^2-1$ & $2a$ \\
\hline $3$ & $4a^3-3a$ & $4a^2 - 1$ \\
\hline $4$ & $8a^4 - 8a^2 + 1$ & $8a^3 -4a$ \\
\hline $5$ & $16a^5 -20 a^3 +5a$ & $16a^4 -12a^2 +1$ \\
\hline $6$ & $32a^6-48a^4+18a^2-1$ & $32a^5-32a^3+6a$ \\
\hline $7$ & $64a^7-112a^5+56a^3-7a$ & $64a^6 -80a^4+24a^2-1$ \\
\hline $8$ & $128a^8-256a^6+160a^4-32a^2+1$ & $128a^7-192a^5+80a^3-8a$\\
\hline
\end{tabular}
\caption{}\label{table}
\end{table}
Note that $X_n(a)$ and $Y_n(a)$ are polynomials in $a$ of degree
$n$ and $n-1$ respectively.
\begin{lem}[J. Robinson's Special Congruence]
\begin{equation}\label{JR1}
Y_n(a)\equiv n\mod (a-1),
\end{equation}
where $a$ and $Y_n(a)$ are as above.
\end{lem}
\begin{proof}
See \cite{matiya}.
\end{proof}
\begin{exam}
By $(\ref{JR1})$ we have that for almost all $a\in \mathbb{Z}$,
$(a-1)\vert(Y_n(a)-n)$. Since $\Z$ is a $D$-ring, then $x-1\vert
Y_n(x)-n$. To have a particular example, take $n=5$, so $Y_5(a) =
16a^4-12a^2+1$, by $(\ref{JR1})$ we have that $a-1\vert
16a^4-12a^2-4$, note that  $x-1 \vert 16x^4-12x^2-4$.
\end{exam}
The following result proved by Julia Robinson, is useful to
show that exponential relations are Diophantine. See \cite{matiya}.
\begin{lem}[J.Robinson]
For all $k\in \mathbb{N}$ we have:
\begin{equation}\label{JR2}
X_n(a)-(a-k)Y_n(a)\equiv k^n\mod{(2ak-k^2-1)}.
\end{equation}
\end{lem}
\begin{exam}
Let $k$ be a non-negative integer. By $(\ref{JR2})$ we have that
for almost all $a \in \mathbb{Z}$, $2ak-k^2-1\vert
X_n(a)-(a-k)Y_n(a)-k^n$, therefore $2xk-k^2-1\vert
X_n(x)-(x-k)Y_n(x)-k^n$. In particular, if $n=7$ then $X_7(a)=
64a^7-112a^5+56a^3-7a$ and $Y_7(a) = 64a^6-80a^4+24a^2-1$. By
$(\ref{JR2})$ we have that
\begin{align*}
2ak-k^2-1\vert
&64a^7-112a^5+56a^3-7a-(a-k)64a^6-80a^4+24a^2-1-k^7\\
               &=-32a^5+32a^3-6a+64a^6k-80a^4k+24a^2k-k-k^7\\
               &=(-1+2ak-k^2)(6a-32a^3+32a^5+k-12a^2k\\
                              &+16a^4k-4ak^2+8a^3k^2-k^3+4a^2k^3+2ak^4+k^5).
\end{align*}
and note that
\[2xk-k^2-1\vert -32x^5+32x^3-6x+64x^6k-80x^4k+24x^2k-k-k^7.\]
\end{exam}
The following Lemma (see \cite{matiya})
provides a relation between the polynomials $X_n(x)$ and $Y_n(x)$.
\begin{lem}
\begin{equation}\label{JR3}
Y_{2n}(a)\equiv 0 \mod {X_n(a)}.
\end{equation}
\end{lem}
\begin{exam}
By $(\ref{JR3})$, for almost all $a \in \mathbb{Z}$ we have that
$X_n(a)\vert Y_{2n}(a)$; and then $X_n(x)\vert Y_{2n}(x)$. If $n=2$, note that
for almost all $a\in \mathbb{Z}$ we have that $2a^2-1 \vert
8a^3-4a$, and $2x^2-1 \vert 8x^3-4x$.
\end{exam}
The following Lemma provides more relations between $X_n(x)$ and
$Y_n(x)$.
\begin{lem}\label{sofisticate}
For $i\geq 1$ we have that:
\begin{align}
Y_{4ni\pm m}(a)&\equiv \pm Y_m(a) \mod{X_n(a)},\\
Y_{4ni+2n\pm m}(a)&\equiv\mp Y_m(a)\mod{X_n(a)}.
\end{align}
\end{lem}
\begin{proof}
See \cite{matiya}.
\end{proof}
\begin{exam}
Let $i\geq 1$, by Lemma \ref{sofisticate} for almost all $a\in \Z$ we have that $X_n(a)| Y_{4ni\pm m}(a) \mp Y_m(a)$, therefore  $X_n(x)|Y_{4ni\pm m}(x)\mp Y_m(x)$.
\end{exam}
\subsection{The ring $\Z[W]$}
We assume the following result from Elementary Number Theory.
\begin{lem}\label{herstein}
Let $p$ be a prime integer and suppose that for some integer $c$
relatively prime to $p$ we can find integers $x$ and $y$ such that
$x^2+y^2 = cp$. Then $p$ can be written as the sum of squares  of
two integers, that is, there exists integers $a$ and $b$ such that
$p=a^2+b^2$.
\end{lem}
\begin{proof}
See \cite[pg~152]{herstein}.
\end{proof}
\begin{thm}[Fermat]\label{fermat}
An odd prime $p$ can be written as $x^2+y^2$ if and only if $p\equiv
1\mod 4$.
\end{thm}
\begin{proof}
See \cite[pg~253]{burton}.
\end{proof}
\begin{exam}
Consider the following set
\[W=\{1/p: p \text{ is prime and } p\equiv 1\mod 4 \text { or } p=2 \}.\]
We take the ring $S=\mathbb{Z}[W]$ and the polynomial $f(x) =
x^2+1$, and we will show that $f(S)\subseteq S^\times$. Let
$\alpha = \frac{a}{b}\in S$. where $a,b \in \mathbb{Z}$ and
$\gcd(a,b)=1$. Note that primes that divide $b$ are primes in
$W$. Note also that the units in $S$ are elements $\frac{c}{d}$ with
$c\equiv 0 \mod p$ and $p\equiv 1\mod 4$. We have that
$f(\alpha)=\frac{a^2+b^2}{b^2}$. Let $p_0$ be a prime such that
$p_0\vert a^2+b^2$, then there exists $c$ such that
$a^2+b^2=cp_0$. By Lemma $\ref{herstein}$, there exist $d$ and
$e$ such that $p_0=d^2+e^2$. By Theorem \ref{fermat}, $p_0\equiv
1\mod 4$. Therefore $f(\alpha)\in S^\times$, this is
$f(S)\subseteq S^\times$. Then, by Proposition \ref{funit} $S$ is
not a $D$-ring. Consequently, $S$ does not satisfy any of the properties $IPP$, $DPP$,
$EPP$ and $SEPP$. Note that $\Z[W]\subseteq \Q$ is a  infinitely generated ring over
$\Z$.
\end{exam}
\subsection{The ring $\Zd$}
Let $d$ be an integer and let $\Zd$ be the subset of complex numbers such that, for every $z\in \Zd$, $z=x+\sqrt{d}y$ with $x,y\in\Z$.
Let $z,w\in \Zd$ and assume $z=x+\sqrt{d}y$ and $w=u+\sqrt{d}v$, we can define arithmetic operations over $\Zd$ as follows:
\begin{align*}
&z+w=(x+u)+\sqrt{d}(y+v),\\
&zw=(xu+dyv)+\sqrt{d}(xv+uy).
\end{align*}
It is easy to see that $\Zd$ with those operations is a domain.
\begin{exam}
If $d=-1$, the domain $\Zd$ is the ring of \textit{Gaussian Integers} $\Z[i]$. If $d=2$, we obtain the domain $\Z[\sqrt{2}]$.
Note that $\Z[i]$ is an Euclidian Domain, therefore it is a $UFD$ with a finite number of units, it is also an infinite domain. By Proposition \ref{infiniteprimes2}, it has an infinite number a prime elements. The ring $\Z[\sqrt{2}]$ is not a $UFD$, because there exist prime elements which are not irreducible elements. Moreover, this ring has an infinite number of units. To see this, note that the equation $x^2-2y^2=1$ has an infinite number of solutions $(x,y)$ because it is a Pell equation. Therefore, the units of $\Z[\sqrt{2}]$ are the element $x+\sqrt{d}y$ such that $x^2-dy^2=1$. Note that $x^2-dy^2=(x+\sqrt{2}y)(x-\sqrt{2}y)$. This example motivates the following definition.
\end{exam}
\begin{defi}
For all $z=x+\sqrt{d}y\in \Zd$ we define the \textit{conjugate} of $z$ as the complex number $\overline{z}=x-\sqrt{d}y$.
\end{defi}
Note that $z=x+\sqrt{d}y\in \Zd$ is a unit if and only if $z\overline{z}=1$. This is: $z$ is a unit in $\Zd$ if and only if $(x,y)$ is solution
of the Pell's equation $x^2-dy^2=1$. Therefore, if $d\geq 2$, the domain $\Zd$ has an infinitely many units. However, it is not known in general for what values of $d$  $\Zd$ is a $UFD$ or not.
The following Lemma shows some elementary properties about the conjugate number.
\begin{lem}\label{lemmaconj}
Let $z,w\in \Zd$. Then:
\begin{enumerate}
\item[(1)] $z\overline{z}\in \Z$;
\item[(2)] $z\in \Z$ if and only if $\overline{z}=z$;
\item[(3)] $\overline{zw}=\overline{z}\cdot\overline{w}$ and $\overline{z+w}=\overline{z}+\overline{w}$;
\item[(4)] $z\overline{w}+\overline{z}w \in \Z$.
\end{enumerate}
\end{lem}
\begin{defi}
Let $f(x)=a_nx^n+\ldots+a_1x+a_0$ with $a_0,a_1,\ldots,a_n\in \Zd$. The \textit{conjugate polynomial} $\C(f(x))$ of $f(x)$ is the polynomial $\C(f(x))=\overline{a_n}x^n+\ldots+\overline{a_1}x+\overline{a_0}$.
\end{defi}
\begin{exam}
Let $f(x)=(1-i)x^2+3ix+1$ in $\Z[i][x]$, then $\C(f(x))=(1+i)x^2-3ix+1$. Let $f(x)=(1-\sqrt{2})x^2-5x+(4-3\sqrt{2})$ in $\Z[\sqrt{2}][x]$, then $\C(f(x))=(1+\sqrt{2})x^2-5x+(4+3\sqrt{2})$.
\end{exam}
Note that every polynomial $f(x)\in \Zd[x]$ can be written as $f(x)=f_1(x)+\sqrt{d}f_2(x)$, where $f_1(x),f_2(x)\in\Z[x]$. Then $\C(f(x))=f_1(x)-\sqrt{d}f_2(x)$.
We have also that if $z\in \Zd$, $\C(z)=\overline{z}$; and for every polynomial $f(x)$ with integer coefficients, $\C(f(x))=f(x)$.
Conversely, if $\C(f(x))=f(x)$ then $f(x)$ is a polynomial with integer coefficients.
The following Proposition shows some elementary properties about the conjugate polynomial.
\begin{prop}\label{propconj}
Let $f(x),g(x) \in \Zd[x]$ and $b\in \Z$. Then:
\begin{enumerate}
\item[(1)] $\C(f(x)+g(x))=\C(f(x))+\C(g(x))$;
\item[(2)] $\C(f(x)g(x))=\C(f(x))\C(g(x))$;
\item[(3)] $\C(f(b))=\overline{f(b)}$;
\item[(4)] $f(x)\C(f(x))\in \Z[x]$;
\item[(5)] $f(x)\C(g(x))+ g(x)\C(f(x))\in \Z[x]$.
\end{enumerate}
\end{prop}
\begin{defi}
Let $f(x)\in \Zd[x]$, we define the \textit{polynomial norm} of $f(x)$ as the polynomial $\Norm(f(x))=f(x)\C(f(x))$. Note that $\deg \Norm(f(x))= 2\deg f(x)$.
\end{defi}
\begin{exam}
Let $f(x)=(1-i)x^2+3ix+1$ in $\Z[i][x]$, then $\Norm(f(x))=[(1-i)x^2+3ix+1][(1+i)x^2-3ix+1]=2x^4-6x^3+11x^2+1$.
Let $g(x)=(1-\sqrt{2})x^2-5x+(4-3\sqrt{2})$ in $\Z[\sqrt{2}][x]$, then $\Norm(g(x))=[(1-\sqrt{2})x^2-5x+(4-3\sqrt{2})][(1+\sqrt{2})x^2-5x+(4+3\sqrt{2})]=-x^4-10x^3+21x^2-40x-2$.
\end{exam}
Note that in the last example, the polynomials $\Norm(f(x))$ and $\Norm(g(x))$ are polynomials with integer coefficients only. This motivates the following result.
\begin{lem}
Let $f(x)\in \Zd[x]$. Then:
\begin{enumerate}
\item[(1)] $\Norm(f(x))=0$ if and only if $f(x)=0$;
\item[(2)] $\Norm(f(x))\in \Z[x]$;
\item[(3)] $\Norm(f(x)g(x))=\Norm(f(x))\Norm(g(x))$;
\item[(4)] for every $a\in \Z$, $\Norm(f(a))=f(a)\overline{f(a)}$.
\end{enumerate}
\end{lem}
\begin{proof}
Immediate from Lemma \ref{propconj}.
\end{proof}
It is already proved in \cite{narki} and \cite{hiroshi} than the domain $\Zd$ is a $D$-ring for every $d\in \Z$. But those proofs are a little complicated and hard to understand. Here, we use the results we have obtained and the above discussion to give an elementary proof that $\Zd$ satisfies $DPP$, consequently $\Zd$ is a $D$-ring for every $d\in \Z$.
\begin{prop}\label{ZdDpp}
For every $d\in \Z$, the ring $\Zd$ satisfies $DPP$. Therefore $\Zd$ is also a $D$-ring.
\end{prop}
\begin{proof}
Let $f(x),g(x)\in \Zd[x]$ be such that for all $k\in \Zd$ $(g(k)\not=0 \Rightarrow g(k)|f(k))$. Consider the polynomials with integer coefficients $F(x)=\Norm(f(x))$ and $G(x)=\Norm(g(x))$. Let $b\in \Z$ such that $G(b)\not=0$ then $g(b)\not=0$. By our choice of $g(x)$, we have that $g(b)|f(b)$ and $\overline{g(b)}|\overline{f(b)}$. By divisibility properties, $g(b)\overline{g(b)}|f(b)\overline{f(b)}$. This implies that $G(b)|F(b)$. We had proven that for every $b\in \Z$, $(G(b)\not=0 \Rightarrow G(b)|F(b))$. Since $\Z$ satisfies $DPP$, $\deg G(x)\leq \deg F(x)$ or $F(x)=0$. Hence $\deg g(x)\leq \deg f(x)$ or $f(x)=0$. In other words, $\Zd$ satisfies $DPP$.
\end{proof}
\begin{cor}
For every $d\in \Z$, the ring $\Z\left[\frac{1+\sqrt{d}}{2}\right]$ satisfies $DPP$. Therefore $\Z\left[\frac{1+\sqrt{d}}{2}\right]$ is a $D$-ring.
\end{cor}
\begin{proof}
Immediately from Proposition \ref{ZdDpp} and Corollary \ref{dSDpp}.
\end{proof}
Note that the argument used to prove that $\Zd$ satisfies $DPP$ is also useful to prove that $\Zdn$ satisfies $DPP$. Therefore, we have the following Corollary.
\begin{cor}
For every $d_1,\ldots,d_n \in \Z$, the ring $\Zdn$ satisfies $DPP$. Therefore $\Zdn$ is a $D$-ring.
\end{cor}

\end{document}